\documentclass[3p,12pt]{elsarticle}
\usepackage[english]{babel}
\usepackage{dsfont}
\usepackage{amssymb}
\usepackage{amsthm}
\usepackage{amsmath}
\usepackage{lineno}
\usepackage{empheq}
\usepackage{graphicx}
\usepackage{booktabs}
\usepackage{lineno}
\usepackage{algorithm}
\usepackage{multicol}
\usepackage[utf8]{inputenc}

\usepackage{caption}
\usepackage{subcaption}
\begin{document}
\begin{frontmatter}
\title{A three-dimensional Laguerre one-way wave equation solver}
\author{Andrew V. Terekhov}
\ead{andrew.terekhov@mail.ru}
\address{Institute of Computational Mathematics and Mathematical Geophysics,
630090, Novosibirsk, Russia}
%\address{Novosibirsk State University, 630090, Novosibirsk,
%Russia}
\begin{abstract}
A finite difference algorithm based on the integral Laguerre transform in time for solving a three-dimensional one-way wave equation is proposed. This allows achieving high accuracy of calculation results. In contrast to the Fourier method, the approach does not need to solve systems of linear algebraic equations with indefinite matrices. To filter the unstable components of a wave field, Richardson extrapolation or spline approximation can be used. However, these methods impose additional limitations on the integration step in depth. This problem can be solved if the filtering is performed not in the direction of extrapolation of the wave field, but in a horizontal plane. This approach called for fast methods of converting the Laguerre series coefficients into the Fourier series coefficients and vice versa. The high stability of the new algorithm allows calculations with a large depth step without loss of accuracy and, in combination with Marchuk-Strang splitting, this can significantly reduce the calculation time. Computational experiments are performed. The results have shown that this algorithm is highly accurate and efficient in solving the problems of seismic migration.
\end{abstract}
\end{frontmatter}
\section{Introduction}
Consider the one-way wave equation
\begin{equation}
\frac{\partial \hat{u}}{\partial z}=\mathrm{i}\frac{\omega}{c}\sqrt{1-\frac{c^2|\mathbf{k}|^2}{\omega^2}}{\hat{u}},
\label{one-way}
\end{equation}
where $\mathrm{i}=\sqrt{-1}$, $\mathbf{x}=(x,y)$, $\hat{u}\equiv \hat{u}(\mathbf{k},z,\omega)$ is a wave component, $\omega$ is the angular frequency,  $\mathbf{k}=\left(k_x,k_y\right)$ is the wave number vector,  $c=c(\mathbf{x},z)$ is the wave velocity, and the vertical direction  $z$ is the direction of extrapolation. Mathematical models based on this pseudodifferential equation have been successfully used to solve problems in seismic prospecting \cite{Claerbout:1985,Symes2008,Biondi2006,Angus2013}, ocean acoustics \cite{Lee2000}, simulation of electromagnetic waves \cite{Leontovich1946,Levy2000}, as well as for setting non-reflecting boundary conditions   \cite{Lindman1975,Engquist1977}. Extrapolation of wave fields in three-dimensional inhomogeneous media is difficult, since it requires simultaneously providing stability, high accuracy, and computational efficiency of the method.

In practice, numerical-analytical algorithms have been widely used (see \cite{Ristow1994,Li1991,Zhang2008,ZHAO2019118,Gazdag1984,Stoffa1990}). Spectral resolution in time and space combined with the discrete  fast Fourier transform (FFT) provides efficient calculations. The accuracy of this class of methods depends on the consistency between the original velocity model and the simplified auxiliary one. High accuracy in simulating wave field kinematics is obtained in the so-called non-stationary phase shift method \cite{Margrave1999}, whereas a more complicated mathematical model has been developed to correctly account for the amplitudes \cite{Zhang2005}. Both approaches are computationally extensive, especially in 3D geometry.

\begin{figure}[!htb]
\centering
        \includegraphics[width=\textwidth]{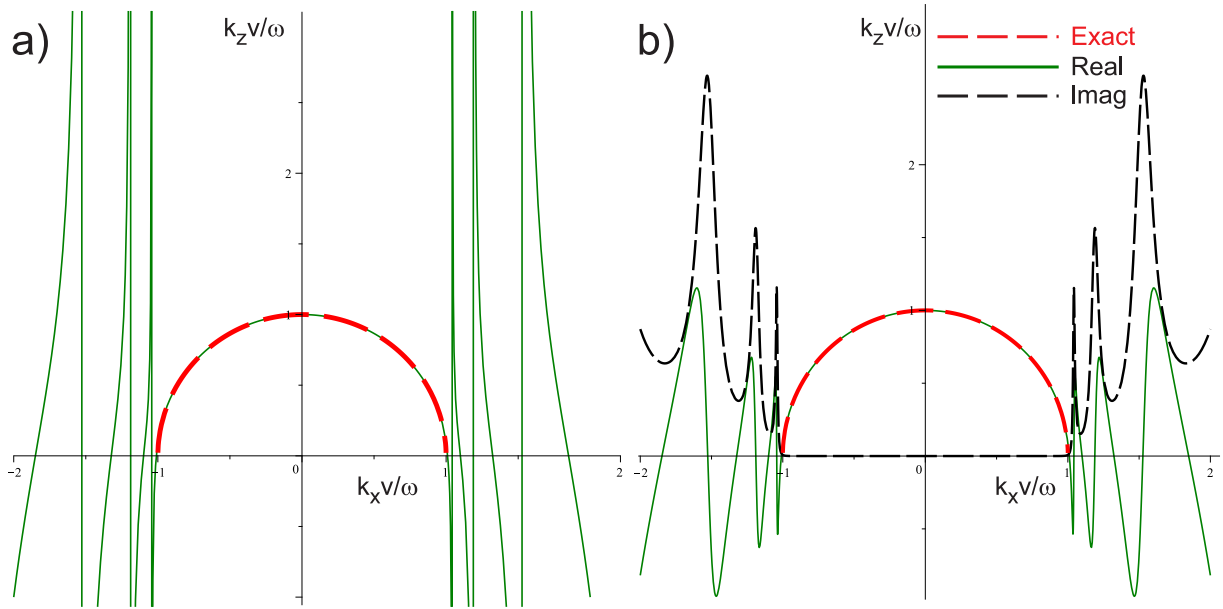}
 \caption{ Padé approximation of the dispersion relation based on (a)~real and (b)~complex coefficients $\gamma_s,\beta_s$ for $n=5$.}
\label{fig:test10}
\end{figure}

Another class of algorithms is based on finite-difference approximations of the spatial derivatives. First, a transition is made from a pseudodifferential equation to a differential equation based on a Padé approximation of the form \cite{Lee01101985}:
\begin{equation}
\sqrt{1-\frac{c^2|\mathbf{k}|^2}{\omega^2}} \approx \mathrm{i}\frac{\omega}{c}\left(1-\sum_{s=1}^n\frac{\beta_sc^2 |\mathbf{k}|^2}{\omega^2-\gamma_sc^2 |\mathbf{k}|^2}\right).
\label{sqr-root-approx}
\end{equation}
 The coefficients $\gamma_s,\beta_s$ , which can be real as well as complex \cite{Lee01101985,Milinazzo1997,Lu1998,Yevick2000}, are determined by minimizing the dispersion error in the least-squares sense. For the real Padé approximation, the denominators in the right-hand side of (\ref{sqr-root-approx}) can take values that are arbitrarily close to zero, which may be a source of solution instability. If  $\gamma_s,\beta_s$ are specified as complex coefficients, the growth of unstable components of the wave field is limited, since it adds a positive imaginary part to the solution of the dispersion equation in the domain of rapidly decaying waves (Fig.~\ref{fig:test10}). However, the strong heterogeneity of the velocity model, as well as the discretization methods used for the original equation, distort the dispersion relation; in some cases this does not allow one to completely remove instability \cite{Amazonas2010}. In addition, the complex Padé approximation is less accurate than the real one \cite{Yevick2000}.

The Marchuk-Strang splitting method is often used to reduce computational work \cite{Marchuk1968,Strang1968}. The problem (\ref{one-way}),(\ref{sqr-root-approx}) transformed from the spectral region to the spatial one takes the following form:
\begin{subequations}
\label{marchuk_decomposition}
  \begin{empheq}[left=\empheqlbrace]{align}
  \label{linz_eq_fourier}
\displaystyle \frac{\partial \tilde{u}_0}{\partial z}=\frac{\mathrm{i}\omega}{c}\tilde{u}_0,\\
  \label{diff_eq_fourier}
\displaystyle \frac{\partial \tilde{u}_{s}}{\partial z}=-\frac{\mathrm{i}\omega}{c}\frac{\beta_sc^2\mathcal{L}}{\omega^2+\gamma_{s}c^2\mathcal{L}}\tilde{u}_{s},\quad s=\overline{1,n},
\end{empheq}
\end{subequations}
where  $\mathcal{L}=\left(\partial^2_x+\partial^2_y\right)$. The sequence of calculations is specified by the initial conditions  $\tilde{u}_0(\mathbf{x},z,\omega)=\tilde{u}(\mathbf{x},z,\omega)$, $\tilde{u}_s(\mathbf{x},z,\omega)=\tilde{u}_{s-1}(\mathbf{x},z+\Delta z,\omega), s=\overline{1,n}$, where the sought-for solution is taken to be $\tilde{u}(\mathbf{x},z+\Delta z,\omega)=\tilde{u}_n(\mathbf{x},z+\Delta z,\omega)$. The finite-difference method for inverting operators of the form $\gamma_{s}c^2\mathcal{L}+\omega^2$ requires significant computational work at the stage of solving the systems of linear algebraic equations. A more efficient approach is based on an alternating direction method \cite{Collino1995,Zhang2008} where special correction procedures are used to decrease the effects of numerical anisotropy \cite{Li1991}.

It was shown in \cite{Costa2013} that a set of complex coefficients $\gamma_s,\beta_s$  for the Padé approximation (\ref{sqr-root-approx}) allows reducing the computational work when using a finite-difference method for equations (\ref{diff_eq_fourier}). This is due to the fact that for high frequencies, $\omega>\omega_0$, the matrix of the finite-difference problem is diagonally dominant. This ensures both high convergence rate of the iterative methods and stability of the direct algorithms used for solving the systems of linear algebraic equations. Unfortunately, for the frequencies $\omega<\omega_0$ the matrices will be ill-conditioned, which significantly decreases the efficiency of the Fourier method for solving three-dimensional problems. To decrease the computational work, it is proposed in papers \cite{Terekhov2017,Terekhov2018}  to use (instead of the Fourier transform in time) the Laguerre transform

\begin{equation}
f(t)=\eta\sum_{m=0}^{\infty}\bar{f}^{m}l_m(\eta t),
\label{series_lag.sum}
\end{equation}
\begin{equation}
\bar{f}^m=\int_{0}^{\infty}f(t)l_m(\eta t) dt,
\label{series_lag.int}
\end{equation}
where $l_m(t)$  are Laguerre functions, and $\eta>0$  is a parameter that controls the convergence rate of the series. In this case, the problem (\ref{one-way}), (\ref{sqr-root-approx}) can be written as a system of equations in the space-time domain \cite{Bamberger1988a}. The functions $l_m(t)$  oscillate on the interval $0<t<4m$ \cite{Temme1990}. Therefore, the Laguerre series with a finite number of terms approximates the function on a finite time interval. To calculate a solution that is non-periodic in time it is not computationally efficient to use Fourier series, since the iterative methods for indefinite systems of linear algebraic equations converge slowly \cite{Ernst2012}. In this case, the Laguerre transform is a good alternative to the Fourier method, since it does not need solving systems of linear algebraic equations with indefinite matrices.

\section{A Laguerre one-way wave equation solver}
\subsection{Discretization of the problem}
Let us transform equation (\ref{marchuk_decomposition}) to the time domain
\begin{subequations}
\label{marchuk_decomposition_time}
  \begin{empheq}[left=\empheqlbrace]{align}
  \label{marchuk_decomposition1}
\displaystyle \frac{\partial u_0}{\partial z}=-\frac{1}{c}\frac{\partial u_0}{\partial t},\\
\label{marchuk_decomposition2}
\displaystyle \left(\frac{1}{c^2}\frac{\partial^2}{\partial t^2}-\gamma_s\mathcal{L}\right)\frac{\partial u_{s}}{\partial z}=\frac{\beta_s}{c}\frac{\partial }{\partial t}\mathcal{L}u_s,\; s=\overline{1,n}.
\end{empheq}
\end{subequations}
Taking into account $f(0)=f'(t)=0$ and assuming \mbox{$\lim_{t\rightarrow \infty}f(t)=\lim_{t\rightarrow \infty}f'(t)=0$}, one can show \cite{Integral_Transform} that
\begin{equation}
\int_{0}^{\infty}\frac{d^{k}f(t)}{d t^k}l_m(\eta t)dt=\left(\frac{\eta}{2}\right)^{k}\bar{f}^m+\Phi_k(\bar{f}^m),
\label{partial_t_lag}
\end{equation}
where for $k=1,2$ we have
\begin{equation}
\Phi_1(\bar{f}^m)\equiv \eta\sum_{k=0}^{m-1}\bar{f}^k, \quad\Phi_2(\bar{f}^m)\equiv\eta^2\sum_{k=0}^{m-1}(m-k)\bar{f}^k.
\label{series_prop}
\end{equation}
 A solution to the one-way wave equation will be sought for in the form of a series of Laguerre functions:
\begin{equation}
u(\mathbf{x},z,t)\approx\sum_{m=0}^{M-1}\bar{u}^m(\mathbf{x},z)l_m(\eta t),
\label{solution_approx}
\end{equation}
where the number of expansion coefficients is chosen according to the required accuracy of approximation of the solution on an interval, $t\in [0,L]$. Applying the Laguerre transform to the system (\ref{marchuk_decomposition_time}), taking into account the properties (\ref{partial_t_lag}), we obtain
\begin{subequations}
\label{marchuk_lag}
  \begin{empheq}[left={ }\empheqlbrace]{align}
\label{linz_eq2}
\displaystyle \frac{\partial \bar{U}^m_0}{\partial z}=-\frac{1}{c}\left(\frac{\eta}{2}\bar{U}_0^m+\Phi_1(\bar{U}^m_0)\right),\\\displaystyle \frac{1}{c^2}\left(\frac{\eta^2}{4}\frac{\partial \bar{U}^m_{s}}{\partial z}+\Phi_2\left(\frac{ \partial \bar{U}^m_{s}}{\partial z}\right)\right)-\gamma_s\mathcal{L}\frac{\partial \bar{U}^m_{s}}{\partial z}\label{diffraction_eq} \\ \nonumber \displaystyle =\frac{\beta_s}{c}\mathcal{L}\left(\frac{\eta }{2} \bar{U}^m_{s}+\Phi_1\left(\bar{U}^m_{s}\right)\right),\  s=\overline{1,n},\;\; m=\overline{0,M-1}.
\end{empheq}
\end{subequations}
An analytical solution for (\ref{linz_eq2}) is (see \cite{Terekhov2019,Sumita1988})
\begin{equation}
\begin{array}{ll}
\displaystyle \bar{U}_0^m(\mathbf{x},z+\Delta z)=\sum_{j=0}^m\left[\bar{U}^{m-j}_0(\mathbf{x},z)-\bar{U}^{m-j-1}_0(\mathbf{x},z)\right]l_{j}(\eta\Delta z/c).
\label{analytic_lag_linz}
\end{array}
\end{equation}
It can be calculated in $O(M\log M)$  arithmetic operations if an FFT algorithm is used to calculate linear convolution \cite{Nussbaumer1982}. The solution to equation (\ref{linz_eq_fourier}) has the form
\begin{equation}
\label{analytic_fourier_linz}
\tilde{u}_0(\mathbf{x},z+\Delta z,\omega)=\exp{\left(\mathrm{i}\omega \Delta z/c\right)}\tilde{u}_0(\mathbf{x},z,\omega),
\end{equation}
and can be calculated in less arithmetic operations than (\ref{analytic_lag_linz}).

Consider a finite-difference approximation for equation (\ref{diffraction_eq}) on a uniform grid with steps $\Delta x,\Delta y,\Delta z$  and nodes  $x_i=(i-1)\Delta x$, $y_j=(j-1)\Delta y$, $z_k=(k-1)\Delta z$, where  $i=\overline{1,N_x}$, $j=\overline{1,N_y}$, $k=\overline{1,N_z}$. The operator $\mathcal{L}$  is approximated  by using the dispersion-relation-preserving finite difference scheme \cite{Tam1993262,Zhang2013511}
\begin{equation}
\begin{array}{ll}
\displaystyle
\mathcal{L}^hf_{i,j,k}=&\displaystyle\ \ \ \left[a_0f_{i,j,k}+\sum_{p=1}^{6}a_p\left(f_{i+p,j,k}+f_{i-p,j,k}\right)\ \right]/ {\Delta x^2}\\&\displaystyle +
\left[a_0f_{i,j,k}+\sum_{p=1}^{6}a_p\left(f_{(i,j+p,k}+f_{i,j-p,k}\right)\right]/{\Delta y^2},
\end{array}
\label{approx_L}
\end{equation}
where  $f_{i,j,k}=f(x_i,y_j,z_k)$, and the coefficients $a_0=-3.12513824$, $a_1=1.84108651$, $a_2=-0.35706478$, $a_3= 0.10185626$, $a_4=-0.02924772$, $a_5=0.00696837$, $a_6=-0.00102952$ provide the twelfth order of accuracy. Using the Crank-Nicholson scheme \cite{CrankNicolson2} to approximate the derivative  $\partial_z$, we obtain the following difference problem:
\begin{equation}
\begin{array}{l}
\left(\mathcal{L}^{h}-\frac{\eta}{(2c_{i,j,k}\gamma_s+\hat{\eta}\Delta z\beta_s)}\right)\bar{U}^{m,s}_{i,j,k+1}=\\\\
\displaystyle  -\frac{1}{{2c_{i,j,k}\gamma_s+{\eta \beta_s}/{2}}}\left[\left(\frac{\eta \beta_s\Delta z}{2}-2c_{i,j,k}\gamma_s\right)\mathcal{L}^{h}\bar{U}^{m,s}_{i,j,k}+
\Delta z\beta_s\mathcal{L}^{h}\Phi_1(\bar{U}^{m,s}_{i,j,k+1}+\bar{U}^{m,s}_{i,j,k})\right.\\\\ \displaystyle \left. +\frac{2}{c_{i,j,k}}\left(\frac{\eta^2}{4}\bar{U}^{m,s}_{i,j,k}-\Phi_2(\bar{U}^{m,s}_{i,j,k+1}-\bar{U}^{m,s}_{i,j,k})\right)\right].
  \end{array}
  \label{lag_matrix}
\end{equation}
The matrix in the left-hand side of the equation is symmetric and definite; therefore, a rapidly converging conjugate gradient method can be used \cite{Golub1989history} to calculate the solution to this system. In the Fourier method, the matrix of the difference problem for equation (\ref{diff_eq_fourier}) will be indefinite, which will require much more calculation time for the iterative algorithms in the Krylov subspace \cite{Ernst2012,VanDerVorst_2003}.

The Laguerre method is based on a real Padé approximation  (\ref{sqr-root-approx}), which may cause instability of the calculations. In the previous versions of the algorithm, the calculations were stabilized by procedures of Richardson extrapolation \cite{Terekhov2017} or spline filtering \cite{Terekhov2018}. However, this approach imposes additional restrictions on the integration step in depth, $\Delta z$. This problem can be solved if the filtering is performed not in the direction of extrapolation of the wave field, but in an  $xy$-plane in the space $(\omega-\mathbf{k})$. For this, it will be necessary to develop procedures for converting the Laguerre series coefficients into the Fourier series coefficients and vice versa.
\subsection{Fast Fourier-Laguerre transform}
To expand functions in the Laguerre series (\ref{series_lag.sum}) it is required to calculate integrals of the form (\ref{series_lag.int}) for rapidly oscillating functions. Since the initial data in the seismic problems are specified with a constant discretization step, no high-accuracy Gaussian quadratures with unevenly distributed nodes can be used \cite{Davis2007}. A new approach for calculating the transformation (\ref{series_lag.int}) was proposed in \cite{Terekhov2019}; it is based on solving the transport equation by a method of separation of variables with subsequent removal of the fictitious periodicity. As a result, the Laguerre series coefficients are calculated from the Fourier series coefficients with the following transformation:
\begin{equation}
\label{main_matirx_lag}
\begin{array}{l}
\left(\bar{f}^0 , \bar{f}^1..., \bar{f}^{M-1}\right)^{\mathrm{T}}=V\left(\tilde{f}_0 , \tilde{f}_1 ,  ... ,\tilde{f}_{(N-1)}\right)^{\mathrm{T}}, \\\\
V:=\left(\frac{\left({-\eta/2-\mathrm{i}k_{j-1}}\right)^{s-1}}{\left({\eta/2-\mathrm{i}k_{j-1}}\right)^{s}}\right)_{s=1,j=1}^{M,N},
\end{array}
\end{equation}
where  $k_j=\frac{2\pi}{L}\left(j-\frac{N}{2}\right)$ for the even $N$  and $k_j=\frac{2\pi}{L}\left(j-\frac{N+1}{2}\right)$  for the odd $N$; $\tilde{f}_p=\sum_{j=0}^{N-1}\exp(\frac{-2\pi jp\mathrm{i} }{N})f_j$. This calculation method for any integer $b\geq 0$  adds a fictitious periodicity of the form \mbox{$f(t)=f(t+bL)$},  which (as shown in \cite{Terekhov2019}) can be removed by using the operation  of double conjugation $\mathbb{Q}^2\left\{\bar{f}^j;L\right\}\equiv \mathbb{Q}\left\{\mathbb{Q}\left\{\bar{f}^j;L\right\};L\right\}$, where
\begin{equation}
\label{conjg_operation}
\mathbb{Q}\left\{\bar{f}^j;\tau\right\}=\sum_{m=0}^{\infty}\left(\bar{f}^m-\bar{f}^{m-1}\right)l_{m+j}(\eta \tau), \ \bar{f}^{-1}\equiv0, \quad j=0,1,....
\end{equation}
The linear correlation for two sequences  of finite length can be efficiently calculated by the FFT algorithm \cite{Nussbaumer1982}.

Let us propose an alternative way of removing the fictitious periodicity.  Based on formula (\ref{analytic_lag_linz}), we introduce into consideration the following "right"{} shift operation of ~$\tau\geq 0$:
\begin{equation}
\label{shift_operation}
\mathbb{S}\left\{\bar{f}^m;\tau\right\}=\sum_{j=0}^m\left(\bar{f}^{m-j}-\bar{f}^{m-j-1}\right)l_{j}(\eta \tau), \ \bar{f}^{-1}\equiv0, \quad m=0,1,....
\end{equation}
If the auxiliary values of $\bar{g}^m$  are calculated by formula (\ref{main_matirx_lag}), the correction
\begin{equation}
\bar{f}^m=\bar{g}^m-\mathbb{S}\left\{\bar{g}^m;L\right\},\quad  m=0,1,...
\label{exclude_period}
\end{equation}
can remove the fictitious periodicity. This is due to the fact that the above difference will take zero values with good accuracy everywhere except for the interval $t \in [0,L]$. This approach is more efficient than the operation $\mathbb{Q}^2$, since only one convolution (\ref{shift_operation}) (instead of two correlations (\ref{conjg_operation})) is to be calculated.

If it is known a priori that a Laguerre series takes zero values at  $t=0$ and $t\geq L$, the inverse transform for $t\in[0,L]$ can be calculated by using, instead of the Laguerre series coefficients, the Fourier series coefficients:

 \begin{equation}
\label{main_matirx_lag_inv}
  \left(\tilde{f}_0 , \tilde{f}_1 ,  ... ,\tilde{f}_{(N-1)}\right)^{\mathrm{T}}=\frac{1}{L}V^{*}
  \left(\bar{f}^0 , \bar{f}^1..., \bar{f}^{M-1}\right)^{\mathrm{T}},
\end{equation}
and then the needed values in the time domain can be calculated by using the inverse discrete Fourier transform.  To make the Laguerre series on the interval $t\in [L,\infty)$   zero with a sufficient accuracy,  the operation $\mathbb{Q}^2\left\{\bar{f}^m;L\right\}$  is applied to its coefficients. The correction of the form (\ref{exclude_period}) cannot be used in this case, since the function approximated by the series is not periodic in the general case.

\begin{figure}[!htb]
\centering
           \includegraphics[width=\textwidth]{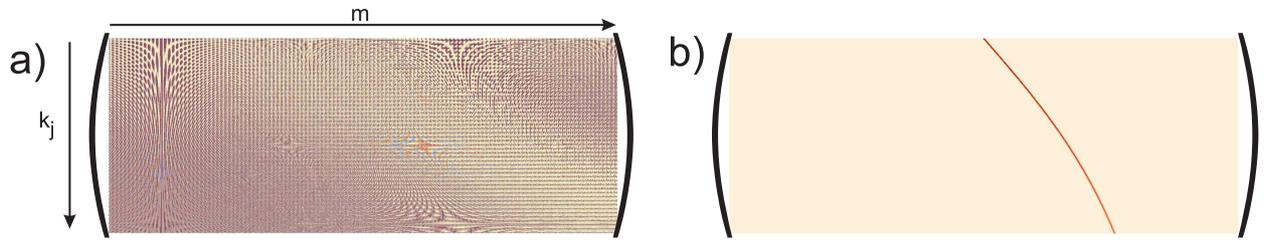}
          \caption{ (a)~Matrix  of the Fourier-Laguerre transform and (b)~the corresponding compressed matrix obtained by the extra-component method.}
\label{Compress1}
\end{figure}

To rapidly calculate the matrix-vector products (\ref{main_matirx_lag}) and (\ref{main_matirx_lag_inv}), an extra-component method \cite{Terekhov2020} can be used to decrease the number of operations from $O(NM)$  to $O(M\log N)$. The main idea of the method is that first, at the preliminary stage, the matrix of the Fourier-Laguerre transform (Fig.~\ref{Compress1}a) is reduced to compact form (Fig.~\ref{Compress1}b)) by using a compression procedure in $O(NM\log M)$  arithmetic operations. Then all multiplications of the vectors by the thus obtained compressed matrix can be performed only in $O(M\log M+24N)$  operations. The large number of matrix-vector products to be calculated allows neglecting the costs of the preparatory stage of the extra-component method.

\subsection{Extrapolation procedure with filtration}
The possibility to rapidly calculate the Laguerre series coefficients into the Fourier series coefficients and vice versa allows using spectral filters to suppress the unstable components of a wave field caused by the real Padé approximation. The above algorithms for converting the Laguerre series coefficients into the Fourier series coefficients and vice versa make it possible to filter the wave field in an  $xy$-plane and perform stable calculations for larger values of $\Delta z$. Let us formulate an elementary step of the wave field extrapolation algorithm.
\begin{algorithm}[H]
  \caption{To calculate the expansion coefficients $\bar{u}^{m}_{i,j,{k+1}}$  of the series (\ref{solution_approx}) with~$\bar{u}^{m}_{i,j,k}$  known at the previous depth layer:}
  \label{alg3}
  \begin{enumerate}
  \item Calculate the grid functions $\bar{U}^{m,s}_{i,j,k+1},$  $m=\overline{ 0,M-1}$     by equation (\ref{lag_matrix})  in the order given by the initial conditions $\bar{U}^{m,1}_{i,j,k}=\bar{u}^{m}_{i,j,k},$ $\bar{U}^{m,s}_{i,j,k}=\bar{U}^{m,s-1}_{i,j,k+1},\; s=\overline{2,n}$ .
    \item Remove the fictitious periodicity by using the operation  $\bar{U}^{m}_{i,j,k+1}=\mathbb{Q}^2\left\{\bar{U}^{m,n}_{i,j,k+1};L\right\}$.
  \item Calculate the Fourier series coefficients   $\bar{U}^{m}_{i,j,{k+1}}\rightarrow\check{u}^{p}_{i,j,k+1}$, using the relation  (\ref{main_matirx_lag_inv}).
  \item Calculate the phase shift  $\tilde{u}^p_{i,j,k+1}=e^{\left(\mathrm{i}\omega_p \Delta z/c\right)}\check{u}^{p}_{i,j,k+1}$, $p=\overline{0,N-1}$.
  \item Make 2d filtering in the  $xy$-plane  $\dot{\tilde{u}}^{p}_{i,j,k+1}=F\tilde{u}^{p}_{i,j,k+1}$, \mbox{$\ p=\overline{0,N-1}$}.
  \item Calculate the Laguerre series coefficients   $\dot{\tilde{u}}^{p}_{i,j,k+1}\rightarrow\bar{u}^{m}_{i,j,{k+1}}$, using the relation   (\ref{main_matirx_lag}).
  \item Remove the fictitious periodicity for the Laguerre series coefficients  $\bar{u}^m_{i,j,k+1}$ formula  (\ref{exclude_period}).
  \end{enumerate}
  \end{algorithm}
At step $4$, formula (\ref{analytic_fourier_linz}) is used instead of (\ref{analytic_lag_linz}) to decrease the computational work. Also, if the velocity model at the  $k$-th level is homogeneous, that is, $c(\cdot,\cdot,z_k)=\chi= const$, the solution of equation (\ref{one-way}) in the space $(\omega-\mathbf{k})$  is calculated analytically \cite{claerbout1976}:
\begin{equation}\label{analyt_phase}
\tilde{u}^{p}_{i,j,k+1}=\underset{\begin{array}{c}
 \scriptstyle k_{x_i}\Rightarrow x_i,\
 \scriptstyle k_{y_j}\Rightarrow y_j
\end{array}}{\mathrm{IDFT_{2D}}}\left\{\exp\left(\mathrm{i}\frac{\omega\Delta z}{\chi}\sqrt{1-\frac{\chi^2|\mathbf{k}_{i,j}|^2}{\omega^2}\ }\right)\underset{\begin{array}{c}
 \scriptstyle x_i\Rightarrow k_{x_i},\
 \scriptstyle y_j\Rightarrow k_{y_j}
\end{array}}{\mathrm{DFT_{2D}}}\left\{\tilde{u}^{p}_{i,j,k}\right\}\right\}.
\end{equation}
At step $5$, a filter is needed which limits the growth of the unstable wave field components \cite{Bunks1995,Zhou1999}. Let us consider another method based on the analytical solution (\ref{analyt_phase}): Let us specify a set of auxiliary velocities $\chi_l$
$$\min_{\mathbf{x}}c(\mathbf{x},z)=\chi_1<\chi_2<...<\chi_b\leq \max_{\mathbf{x}}c(\mathbf{x},z_k),$$ and use the following filter  $F$:
\begin{equation}
\dot{\tilde{u}}^p_{i,j,k+1}=\underset{\begin{array}{c}
 \scriptstyle k_{x_i}\Rightarrow x_i,\
 \scriptstyle k_{y_j}\Rightarrow y_j
\end{array}}{\mathrm{IDFT_{2D}}}\left\{G(\omega_p,k_{x_i},k_{y_j},\chi_l)\underset{\begin{array}{c}
 \scriptstyle x_i\Rightarrow k_{x_i},\
 \scriptstyle y_j\Rightarrow k_{y_j}
\end{array}}{\mathrm{DFT_{2D}}}\left\{\tilde{u}^{p}_{i,j,k+1}\right\}\right\},\; \chi_l \leq c(x_i,y_j,z_{k}) \leq \chi_{l+1},
\label{filter1}
\end{equation}
where
\begin{equation}
\label{filtration}
G(\omega,k_x,k_y,\chi)=\left\{
\begin{array}{ll}
\exp\left(-\frac{\omega\Delta z}{\chi}\sqrt{\frac{\chi^2|\mathbf{k}|^{2}}{\omega^2}-1}\ \right), &|\mathbf{k}|\geq \frac{\omega}{\chi},\\
1, &|\mathbf{k}|<\frac{\omega}{\chi}.
\end{array}
\right.
\end{equation}
Using a small set of velocities by analogy with a method called PSPI \cite{Gazdag1984} can significantly decrease the filtering time. Although the parasitic waves will be incompletely suppressed, this method will be sufficient for solving practical problems of seismic with high accuracy.

\section{Numerical Experiments}
Consider some computational experiments to assess the accuracy, stability, and efficiency of the above-proposed procedures. The numerical procedures were implemented in \mbox{Fortran-2008} using the Intel MKL library. The calculations were performed on a  supercomputer at Novosibirsk State University. The supercomputer has Intel Xeon Gold 6248 twenty-core processors operating at 2.5~GHz. Each computational node contains four processors and 384~GB of RAM.
\subsection{Fourier-Laguerre transform}
\begin{figure}[!htb]
\centering
\includegraphics[width=\textwidth]{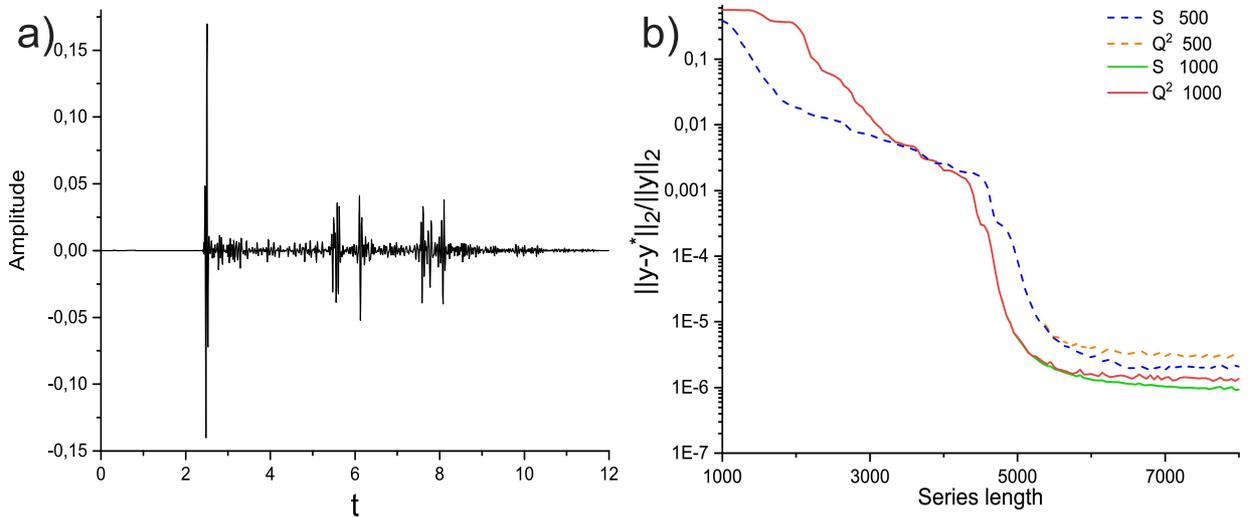}
\caption{a)~First trace from seismograms for the velocity model Sigsbee2A, b)~error of seismic trace approximation by a Laguerre series for various methods of removing the fictitious periodicity for $\eta=500$  and $1000$.}
\label{pic:test3}
\end{figure}

In the algorithm being proposed, a solution to the one-way wave equation is calculated in the form of a Laguerre series. Then the unstable wave field components are filtred in the spectral domain $(\omega-\mathbf{k})$. The conversion of the Laguerre series coefficients into the Fourier series coefficients and vice versa is performed at each depth layer. Therefore, the possible insufficient accuracy of this procedure may cause numerical instability. Let us compare the accuracy of the operation $\mathbb{Q}^2$  and formula  (\ref{exclude_period}) developed to remove the fictitious periodicity emerging when a function is expanded into a Laguerre series based on the transform (\ref{main_matirx_lag}). The first seismic trace (Fig.~\ref{pic:test3}a) from the set of data for the Sigsbee 2A velocity model is taken as a function to be approximated \cite{Paffenholz}. Fig.~\ref{pic:test3}b shows that the accuracy of formula (\ref{exclude_period})  is somewhat higher than that of the operation $\mathbb{Q}^2$. Thus, the new method for removing the fictitious periodicity when using a Laguerre basis instead of a trigonometric basis makes it possible to halve the computational costs of step 7 in Algorithm 1 without loss of accuracy.

\subsection{2d velocity models }
\begin{figure}[!htb]
\centering
           \includegraphics[height=0.89\textheight]{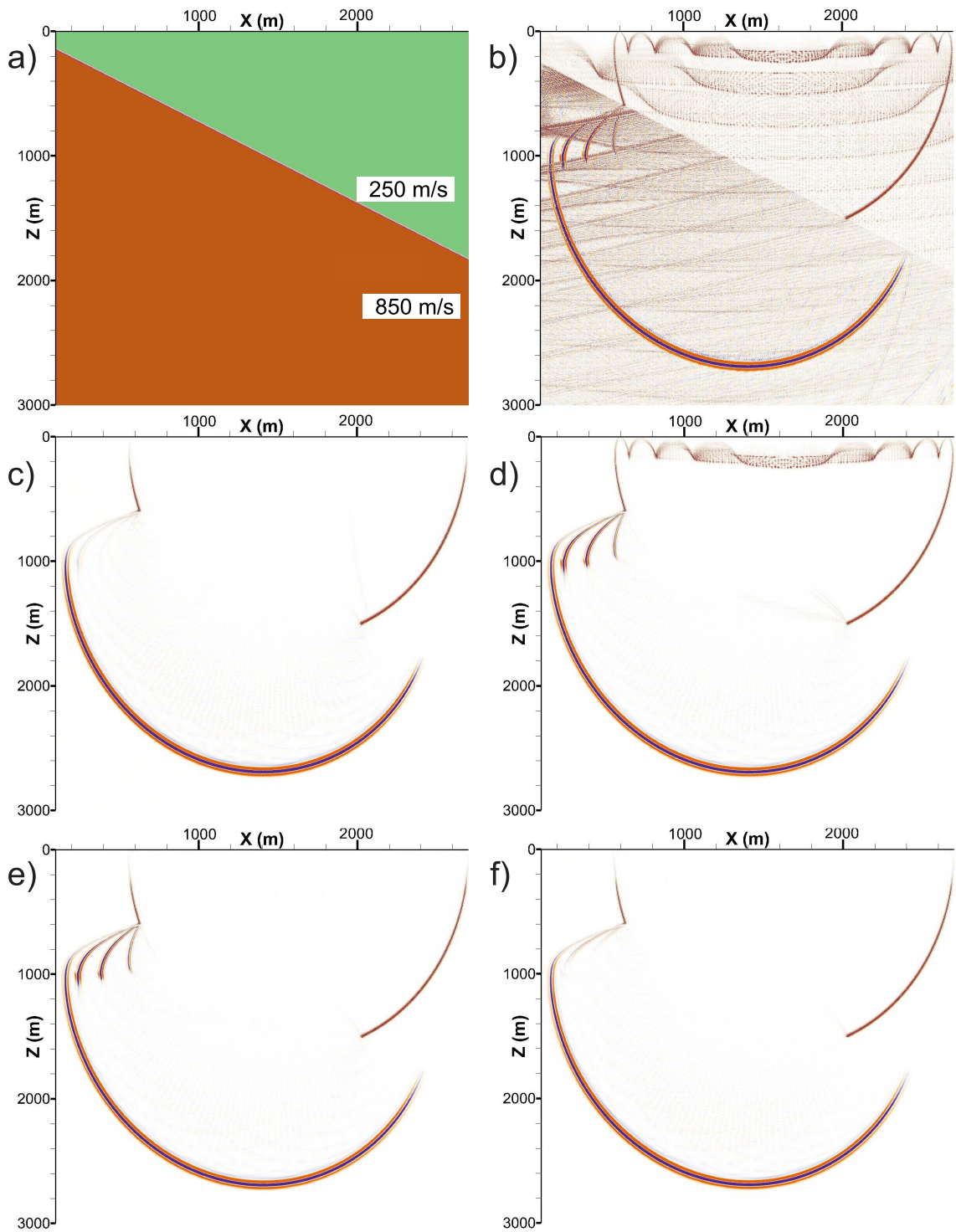}
           \caption{Impulse responses for (a)~incline velocity model. Fourier method with (b)~real and (c)~complex Padé approximations, (d)~Laguerre method without and with filtration for velocities (e)~$\chi_1=250$~m/s, (f)~$\chi_1=250$~m/s, $\chi_2=850$~m/s.}
\label{fig:test1}
\end{figure}

Let us estimate the effect of the filtering procedure implemented in Algorithm~1 on the quality of the solution. For the two-dimensional velocity model (Fig.~\ref{fig:test1}a), we calculate a wave field from a point source whose time dependence is set in the form of a Ricker wavelet with a dominant frequency of  $20$Hz. The number of Laguerre series coefficients $M=2048$  for  $\eta=600$, and for the Fourier series coefficients $N=1024$. Fig.~\ref{fig:test1}b is a snapshot of a wave field calculated by the Fourier method for a real Padé approximation which shows parasitic noise caused by the growth of the unstable solution components. Fig.~\ref{fig:test1}c presents a solution by the Fourier method for complex coefficients that limit the instability.

The solution obtained by the Laguerre method is shown in Fig.~\ref{fig:test1}d. Here the noise level is somewhat lower than that for the Fourier method with real coefficients (Fig.~\ref{fig:test1}b), but much higher than for the Fourier method with complex coefficients (Fig.~\ref{fig:test1}c). Thus, the Laguerre algorithm is conditionally stable, that is,  the unstable components of the wave field can increase, but not exponentially. The filtering procedure for one (Fig.~\ref{fig:test1}e) and two velocities (Fig.~\ref{fig:test1}f) can significantly reduce the numerical noise in the Laguerre method to the level of the Fourier method with complex coefficients. Numerous tests have also confirmed that the repeated conversion of the Laguerre coefficients to the Fourier coefficients and vice versa is a stable procedure.

\begin{figure}[!htb]
\centering
           \includegraphics[width=\textwidth]{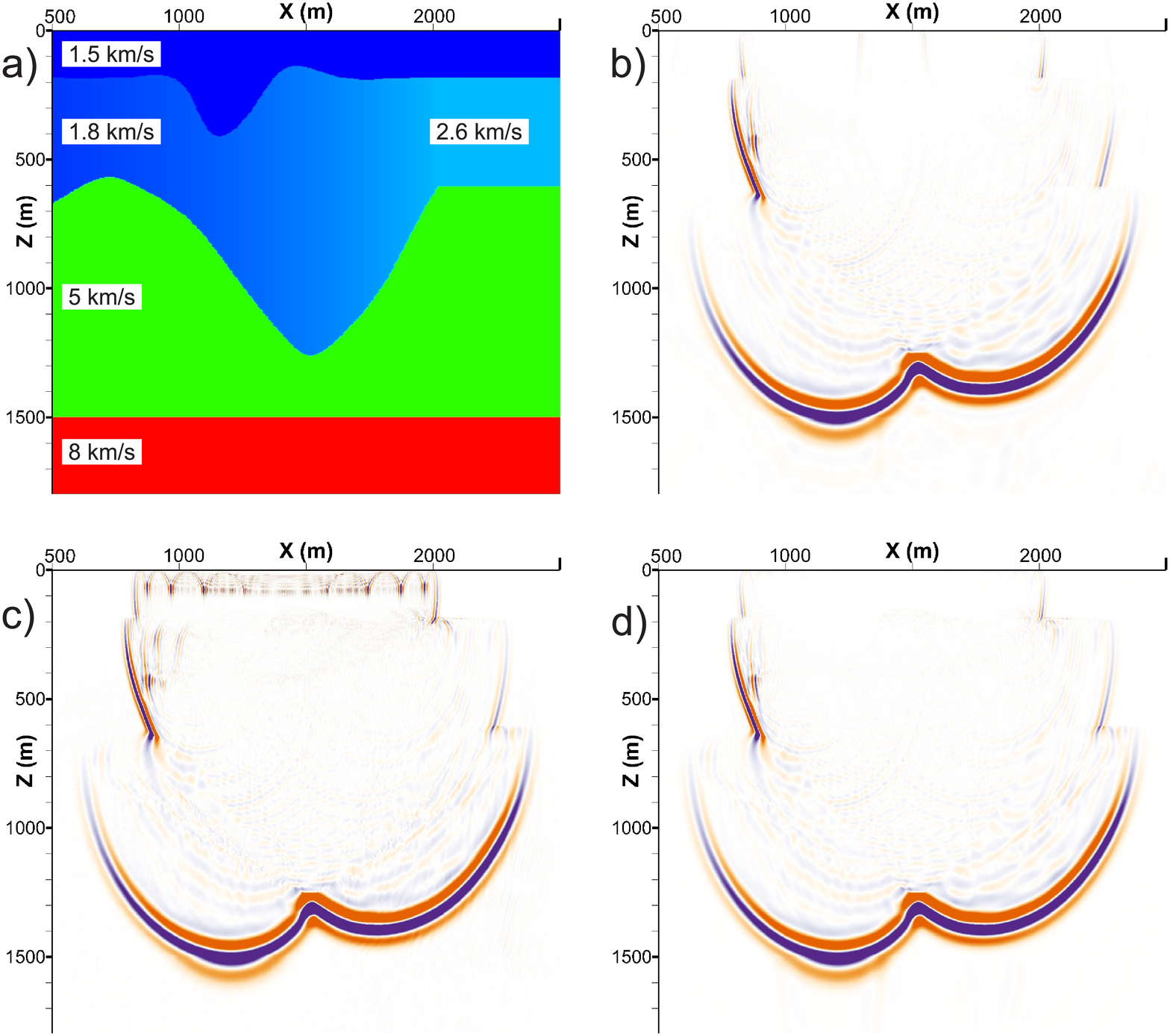}
          \caption{ Impulse responses for (a) syncline velocity model. Fourier method with (b) real and (c) complex Pade-approximations, (d) Laguerre method with filtering for velocities  $\chi_1(z)=\min_{x}c(x,z)$  and $\chi_2(z)=\sum_{i=1}^{N_x}c(x_i,z)/N_x$}
\label{fig:test2}
\end{figure}

Let us consider in Fig.~\ref{fig:test2}a a more complicated velocity model of a medium including a syncline. The numerical calculations were performed on a grid with steps $\Delta x=25$m and  $\Delta z=5$m. The number of the Laguerre series coefficients for $\eta=400$  $M=2000$ , and for the Fourier series $N=1024$. A time dependent point source was specified in the form of a Ricker wavelet with a dominant frequency of $30$Hz. Fig.~\ref{fig:test2}b shows that the solution obtained by the Laguerre method without filtering has much more numerical noise than the solution (Fig.\ref{fig:test2}c) obtained by the Fourier method for a complex Padé approximation. However, the filtering of the unstable components of the wave field for velocities $\chi_1(z)=\min_{x}c(x,z)$   and $\chi_2(z)=\sum_{i=1}^{N_x}c(x_i,z)/N_x$   allows significantly reducing the numerical noise of the Laguerre method (Fig.~\ref{fig:test2}d).

Let us calculate the reflection coefficient for constructing a depth image based on prestack migration \cite{Gray2001}. For each source $s$  with coordinate $\mathbf{x}_s$   the receiver wavefields $R(\mathbf{x}_s,\mathbf{x},z,\omega)$  and the source wavefields $S(\mathbf{x}_s ,\mathbf{x},z,\omega)$  are to be extrapolated from the surface to the depth. After calculating the wave fields, an approximate value of the reflection coefficient can be estimated as follows:
\begin{equation}
\mathrm{I_d}(\mathbf{x},z)=\sum_{\mathbf{x}_s}\sum_{\omega}\frac{R(\mathbf{x}_s,\mathbf{x},z,\omega)
S^{*}(\mathbf{x}_s,\mathbf{x},z,\omega)}{S(\mathbf{x}_s,\mathbf{x},z,\omega)
S^{*}(\mathbf{x}_s,\mathbf{x},z,\omega)+\varepsilon(\mathbf{x}_s,z)},
\label{deconv_imaging}
\end{equation}
where $\varepsilon(\mathbf{x}_s,z)=10^{-3}\max_{\mathbf{x},\omega}|S(\mathbf{x}_s,\mathbf{x},z,\omega)|^2$ is a regularization parameter avoiding the division by zero or near-zero values of the denominator. Let us also consider an image condition that is more stable to noise:
\begin{equation}
\mathrm{I_c}(\mathbf{x},z)=\sum_{\mathbf{x}_s}\sum_{\omega}R(\mathbf{x}_s,\mathbf{x},z,\omega)S^{*}(\mathbf{x}_s,\mathbf{x},z,\omega).
\label{conv_imaging}
\end{equation}
A shortcoming of this condition is that $\mathrm{I_d}$  is nondimensional, which is true from the physical point of view, while the reflection coefficient $\mathrm{I_c}$  has an inappropriate dimension. The quality of depth images obtained by methods of calculating the reflection coefficient is considered in papers \cite{Valenciano2005,Guitton2007,Schleicher2008}.

\begin{figure}[!htb]
\centering
           \includegraphics[width=\textwidth]{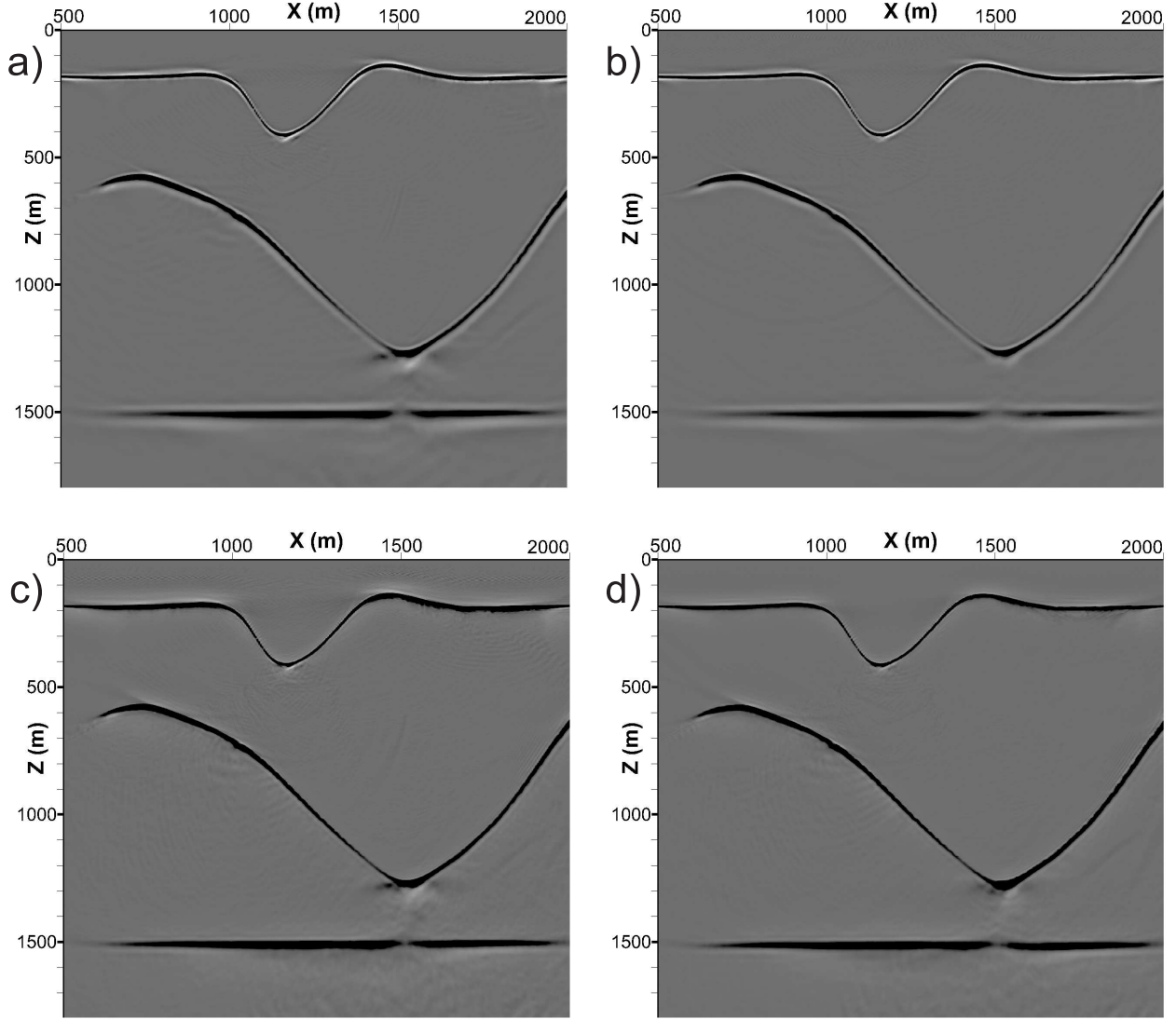}
         \caption{Comparison of the migration results obtained by the Fourier method: (a)~$\mathrm{I_c}$ , (c)~$\mathrm{I_d}$  and the Laguerre method with filtering: (b)~$\mathrm{I_c}$ , (d)~$\mathrm{I_d}$}
\label{fig:test3}
\end{figure}

For the model in Fig.~\ref{fig:test2}a, depth images $\mathrm{I_d}$  and $\mathrm{I_c}$  were calculated for eighty sources. The theoretical seismograms were obtained by using a Gaussian beam algorithm \cite{GaussianBeams_popov,Cerveny1985} implemented in the package Siesmic Unix. Fig.~\ref{fig:test3} shows that the quality of the images calculated by the Laguerre method is comparable to that of the Fourier method. This suggests that the filtering procedure being used provides sufficient suppression of the parasitic harmonics for a significantly non-homogeneous velocity model even when using two filtering velocities.

\begin{figure}[!htb]
\centering
           \includegraphics[width=\textwidth]{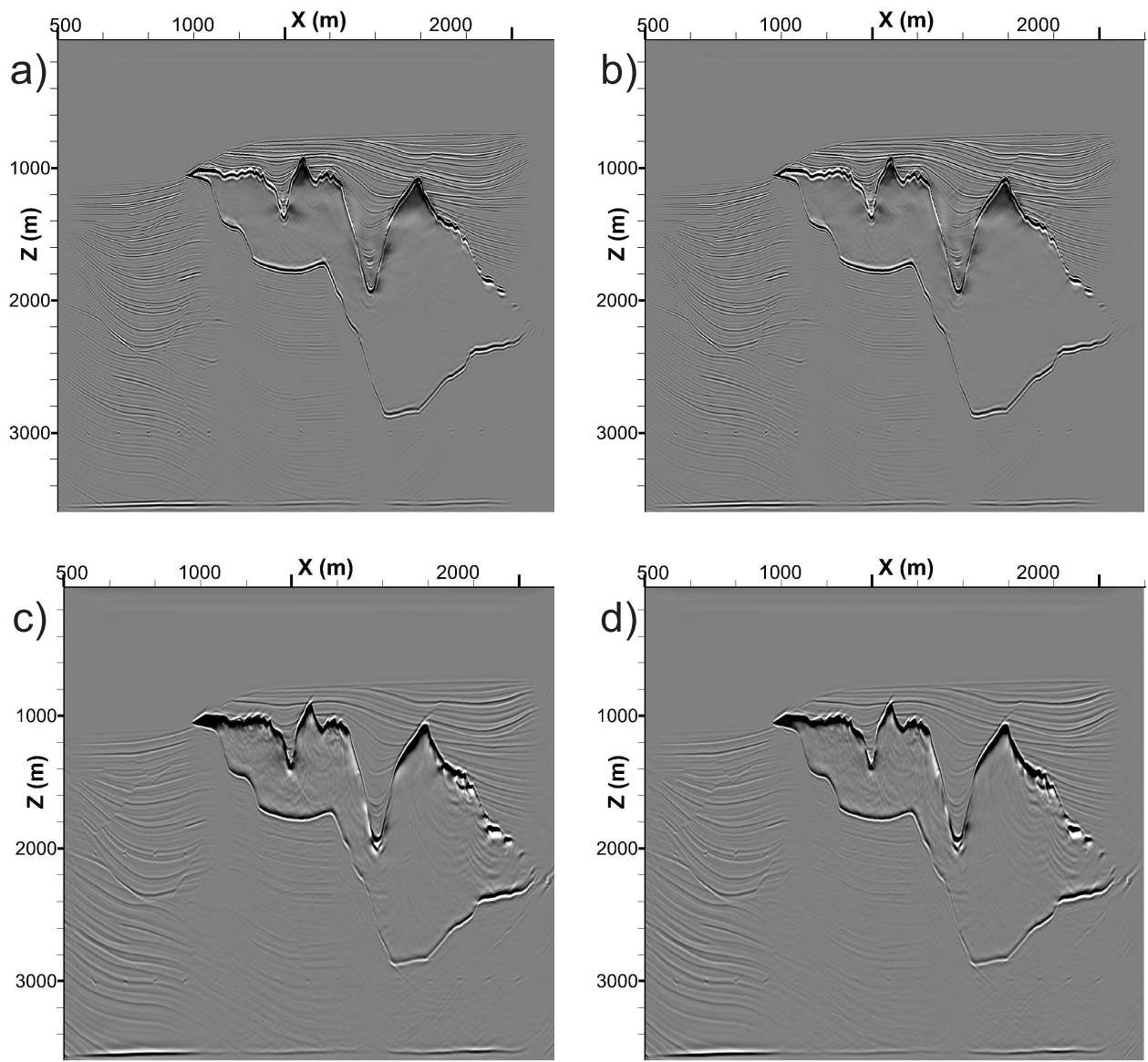}
         \caption{Comparison of migration results obtained by the Fourier method: (a)~$\mathrm{I_c}$, (c)~$\mathrm{I_d}$ and the Laguerre method with filtering: (b)~$\mathrm{I_c}$, (d)~$\mathrm{I_d}$}
\label{fig:test4}
\end{figure}

Consider a solution obtained with a velocity model called Sigsbee2A~\cite{Paffenholz}, for which filtering of the unstable components of the wave field is made with a slightly larger number of auxiliary velocities than in the previous case:
\begin{equation}
\begin{array}{ll}
\chi_1(z)=\min_{\mathbf{x}}c(\mathbf{x},z),& \chi_2(z)=\chi_1(z)+\Delta  \chi(z),\\ \chi_3(z)=\chi_1(z)+2\Delta \chi(z),&  \chi_4(z)=\chi_1(z)+3.2\Delta  \chi(z),\\
\Delta \chi(z)=\left(\max_{\mathbf{x}}c(\mathbf{x},z)-\min_{\mathbf{x}}c(\mathbf{x},z)\right)/3.5.
\end{array}
\label{filt_vel}
\end{equation}
The calculations were performed on a grid with steps  $\Delta x=25$m and  $\Delta z=7.5$m, the number of the Laguerre series coefficients for $\eta=600$ $M=3000$, and for the Fourier series $N=1024$. The depth images for  $500$ sources calculated by the Laguerre method with filtering and the Fourier method coincide (Fig.~\ref{fig:test4}). This makes it possible to conclude that algorithm 1 has sufficient accuracy.

In the two-dimensional case, the Laguerre algorithm has no obvious computational advantages over the Fourier method. This is explained by the fact that in both cases a method based on an LU (or $\mathrm{LL}^{\mathrm{T}}$ ) expansion is used for efficiently solving the systems of linear algebraic equations \cite{Golub1989}. Being a variant of the Gauss method, the LU algorithm calculates the solution in a number of operations that does not depend on the condition number of the matrix. For the Laguerre method, at each depth layer the LU factorization of the matrix (\ref{lag_matrix}) is performed once for all the right-hand sides. The Fourier method for a complex Padé approximation requires no filtering procedure, but the LU factorization is performed not only for each depth layer, but also for each frequency.

In the three-dimensional case, the Laguerre method has a significant advantage in solving the systems of linear algebraic equations, since the matrix is  definite and well-conditioned. Hence, the iterative algorithms will converge much faster than for the Fourier method. In the latter, solving the Helmholtz-type equations is difficult from both theoretical and practical points of view \cite{lahaye_modern_2017}.
\subsection{3D SEG/EAGE Salt  model  }
Let us calculate a wave field from a point source for a 3D velocity model called SEG/EAGE Salt \cite{aminzadeh19973}. The number of the Laguerre and Fourier series coefficients $M=600$   and $N=512$, respectively.  The conversion parameter (\ref{series_lag.int}) was $\eta=500$. The calculations were made on a grid with steps $\Delta x=\Delta y=20$m and $\Delta z=10$m. The total numbers of grid nodes including the  fictitious absorbing Tarper-type boundaries \cite{Cerjan1985} $N_x=N_y=676$ and $N_z=448$ . The point source was specified as a Ricker wavelet with a dominant frequency of  $20$Hz.

In the three-dimensional case, at each depth layer it is necessary to solve 2D difference problems (\ref{lag_matrix}) whose matrix is symmetric and positive definite with a large number of zero elements. To solve such systems of linear algebraic equations, the conjugate gradient method is used. The following stopping criterion for the iterative process was chosen:  $||\mathbf{r}_i||_2/||\mathbf{f}||_2\leq \varepsilon=10^{-6}$, where $\mathbf{r}_i$  is the residual of the  $i$-th iteration of the conjugate gradient method and $\mathbf{f}$   is the vector of the right-hand side of the problem (\ref{lag_matrix}). The use of lower accuracy stopping criteria has caused instability in the wave field extrapolation, which is explained by the loss of consistency of the operators $\mathcal{L}^{h}$  in the right- and left-hand sides of equation (\ref{lag_matrix}).

\begin{figure}[!htb]
\centering
           \includegraphics[width=0.5\textwidth]{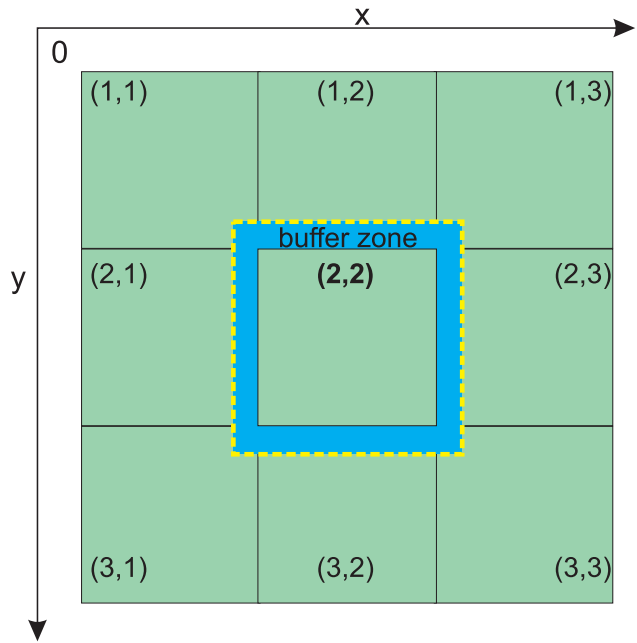}
          \caption{Calculation domain decomposition to decrease the number of iterations of the conjugate gradient method}
\label{DomainDecomposition}
\end{figure}

\begin{figure}[!htb]
\centering
\includegraphics[width=\textwidth]{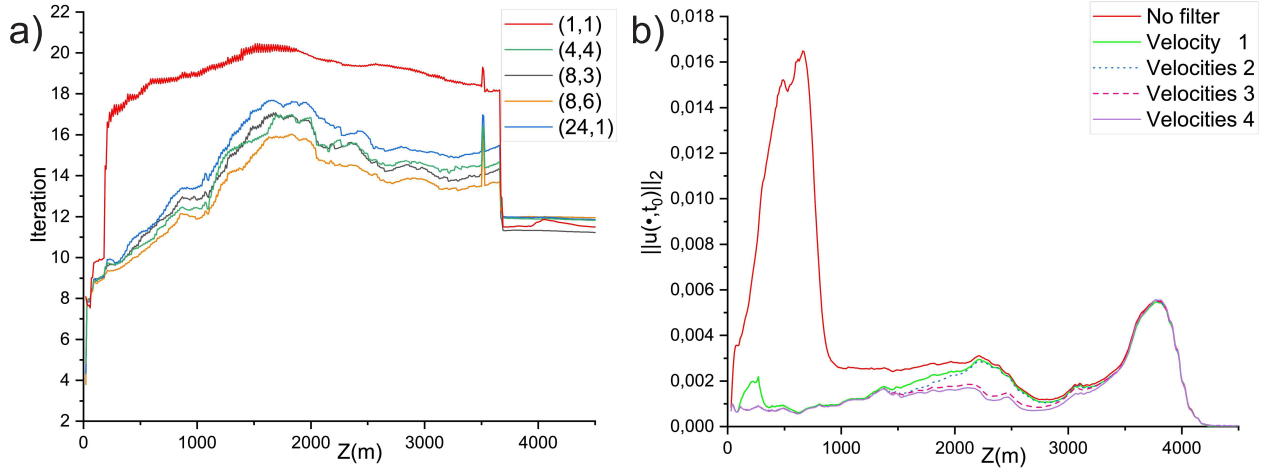}
\caption{For 3D SEG/EAGE Salt velocity model: (a)~the number of iterations of the conjugate gradient method versus the layer depth and various configurations of the domain decomposition method, (b)~wave field energy at  $t=2$~s for the Laguerre method with filtering.}
\label{fig:test7}
\end{figure}

\begin{figure}[!htb]
\centering
           \includegraphics[height=0.89\textheight]{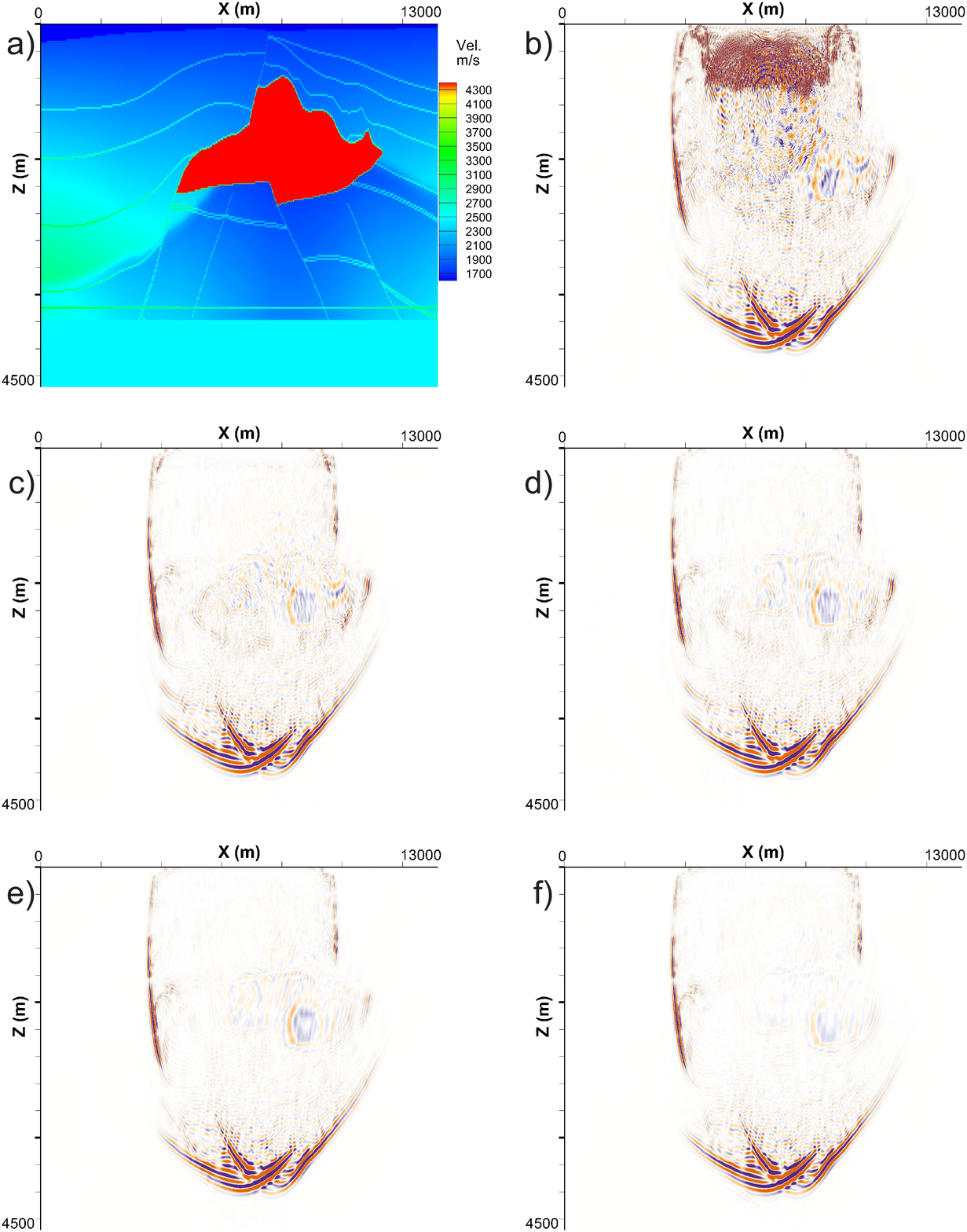}
 \caption{Sections of 3D impulse response in plane $y=6500$~m for (a) SEG/EAGE Salt velocity model. 3D Laguerre method: (b)~without filtering, for (c)~one, (d)~two, (e)~three, and (f)~four filtering velocities in (\ref{filt_vel}).}
\label{fig:test5}
\end{figure}

\begin{figure}[!htb]
\centering
           \includegraphics[height=0.89\textheight]{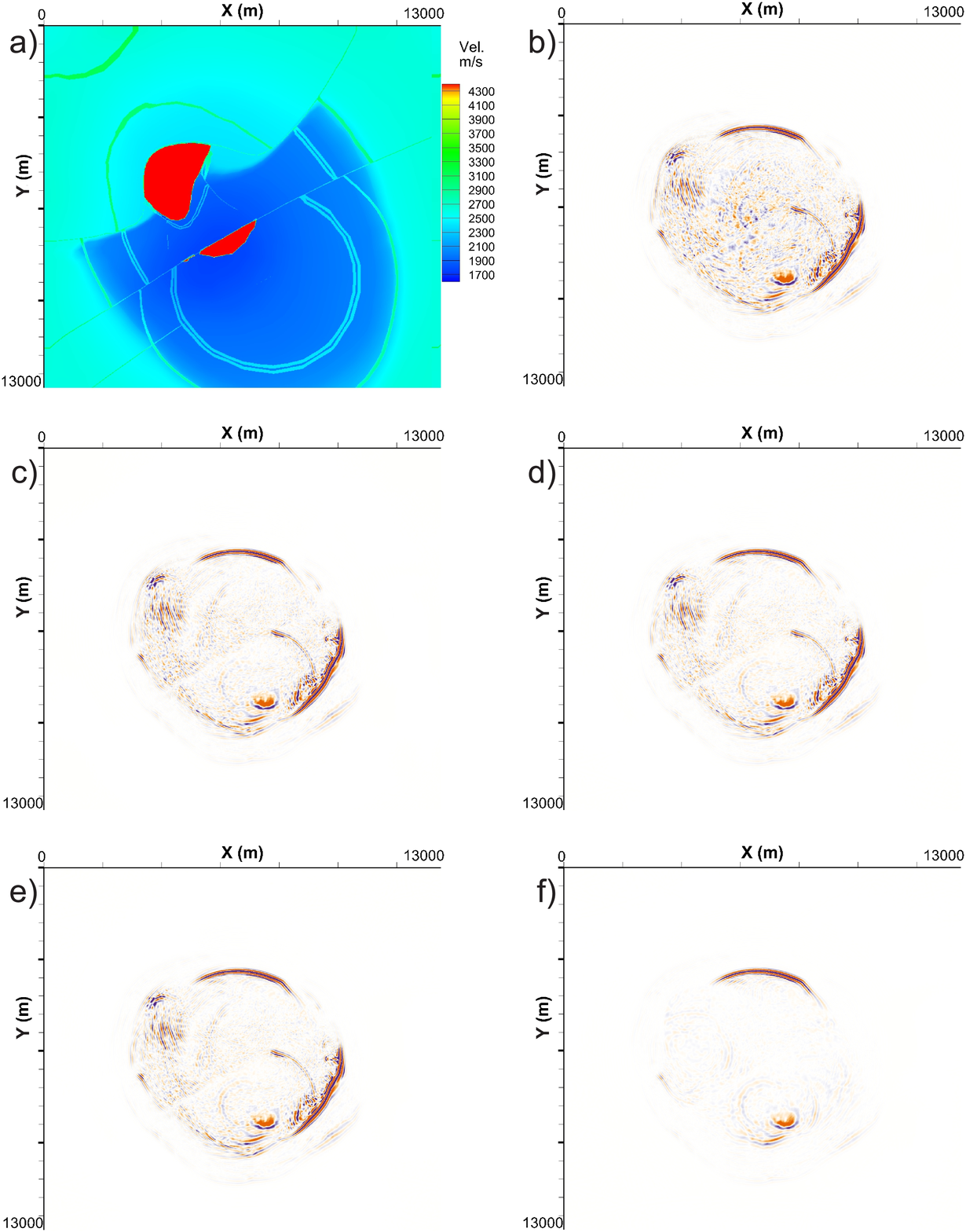}
 \caption{Sections of 3D impulse response in plane $z=1500$~m for (a) SEG/EAGE Salt velocity model. 3D Laguerre method: (b)~without filtering, for (c)~one, (d)~two, (e)~three, and (f)~four filtering velocities in (\ref{filt_vel}).}
\label{fig:test6}
\end{figure}

The use of preconditioning procedures, such as incomplete Cholesky factorization \cite{Golub1989}, makes it possible to decrease the number of iterations by several times, but the total calculation time increases by five to six times due to the additional computational costs. Another approach is based on the fact that the strong diagonal dominance of the matrix of the problem (\ref{lag_matrix}) causes the difference Green function to rapidly tend to zero with distance from the source. This makes it possible to decrease the number of iterations of the conjugate gradient method based on the domain decomposition method as follows \cite{Terekhov2015206}: First the entire 2D domain on the depth layer to be calculated is divided into several non-overlapping rectangular subdomains with possibly equal areas (Fig.~\ref{DomainDecomposition}).
A buffer zone several nodes in width is formed around each subdomain (for the velocity model 3D SEG/EAGE Salt the buffer zone width was eight nodes). In the buffer zone the right-hand side of equation (\ref{lag_matrix}) is set equal to zero, and the values of the right-hand side at the other grid nodes are equal to those at the nodes of the initial domain.  Then the difference boundary value problems with homogeneous Dirichlet conditions are solved independently in each subdomain. Since  the equations are linear, the solution in the entire calculation domain is taken to be equal to the sum of the solutions in the subdomains and the buffer zones. As shown in Fig.~\ref{fig:test7}a, this calculation method makes it possible to reduce the number of iterations by $30-40\%$  with minimal additional work, as well as naturally organize parallel computations. The decrease in the calculations is due to the fact that for many local subproblems the condition number, which depends on $\max_{\mathbf{x}}{c(\mathbf{x},z_k)}/\min_{\mathbf{x}}{c(\mathbf{x},z_k)}$, will be smaller than that for the entire calculation domain. The experiments have shown (see Figs.~\ref{fig:test5} and \ref{fig:test6}) that the calculation method being considered does not introduce considerable amplitude or phase errors. This algorithm cannot be used to solve the difference problem (\ref{diff_eq_fourier}) with the Fourier method, since the corresponding matrix does not have the property of significant diagonal dominance for all given values of~$\omega$. In this case, preconditioned iterative methods in the Krylov subspaces \cite{Golub1989} should be used; their implementation, however, will require computational work that is several orders of magnitude greater.

One can see in Fig.~\ref{fig:test7}a that the Laguerre method needs less iterations to solve the problem (\ref{lag_matrix}) than the Fourier method (see data in \cite{Costa2013}). The total number of iterations depends on the grid step, the velocity model of the medium on the $k$-th layer and does not depend on the number of the expansion coefficient of the Laguerre series (\ref{solution_approx}) being calculated. The calculation time required for filtering and converting the Laguerre series coefficients into the Fourier series coefficients and vice versa constituted $\approx30\%$  of the total calculation time. As in the 2D calculations, Figs.~\ref{fig:test5} and \ref{fig:test6} show that the filtering allows suppressing the energy of parasitic waves (Fig.~\ref{fig:test7}a). However, for the splitting methods with subsequent correction of numerical anisotropy \cite{Zhang2008,Li1991} no auxiliary velocity models can be used to reproduce correctly the dispersion relation for all wave field components in question. There is no such a problem in the approach being proposed.

\section{Conclusions}
A new algorithm for wave field extrapolation down from the surface based on the one-way wave equation was proposed in the present paper. In contrast to the Fourier transform, the integral Laguerre transform in time considerably decreases the computational work in solving the difference subproblems in the 3D case. The above practical calculations have shown that from ten to twenty iterations of the conjugate gradient method must be performed to solve the systems of linear algebraic equations. With a small number of iterations no efficient preconditioning is possible; nevertheless, the computational work can be reduced by a third in the case of a simplified version of the domain decomposition method.
With the Laguerre method, to transform the pseudodifferential equation into the differential one a real Padé approximation (which is unstable) is used. This requires to limit the growth of the unstable wave field components on each depth layer. In the previous works, depth filtering in the spatial domain has been used for this purpose. This imposes additional restrictions on the integration step in depth. In the method proposed, these restrictions are weakened considerably, since the filtering is performed in a horizontal plane with a spectral low-pass filter. Numerous publications have been devoted to the construction and implementation of this class of filters. Therefore, the main difficulty was to develop fast methods for converting the Laguerre series coefficients into the Fourier series coefficients and vice versa. The extra-component algorithm, which was proposed recently  for fast calculation of matrix-vector products, allows solving this problem with high efficiency, since even for complex models of various media the filtering takes no more than thirty percent of the total calculation time. As a result, the total efficiency of the method has increased by at least an order of magnitude. The above practical calculations have shown that the accuracy and stability of the Laguerre algorithm are comparable to those of the Fourier method for a complex Padé approximation, but the former algorithm is much more economical in the number of operations.
\newpage

\newpage
\bibliography{base}

\begin{thebibliography}{10}

\bibitem{Claerbout:1985}
J.~F. Claerbout.
\newblock {\em Imaging the Earth's Interior}.
\newblock Blackwell Scientific Publications, Inc., Cambridge, MA, USA, 1985.

\bibitem{Symes2008}
P.~Shen and W.~Symes.
\newblock Automatic velocity analysis via shot profile migration.
\newblock {\em Geophysics}, 73(5):VE49--VE59, 2008.

\bibitem{Biondi2006}
B.~Biondi.
\newblock {\em 3D Seismic Imaging}.
\newblock Society of Exploration Geophysicists, 2006.

\bibitem{Angus2013}
D.~A. Angus.
\newblock The one-way wave equation: A full-waveform tool for modeling seismic
  body wave phenomena.
\newblock {\em Surveys in Geophysics}, 35(2):359--393, 2013.

\bibitem{Lee2000}
D.~Lee, A.~D. Pierce, and Er-C. Shang.
\newblock Parabolic equation development in the twentieth century.
\newblock {\em Journal of Computational Acoustics}, 08(04):527--637, 2000.

\bibitem{Leontovich1946}
M.~Leontovich and V.~Fock.
\newblock Solution of the problem of propagation of electromagnetic waves along
  the earth’s surface by the method of parabolic equation.
\newblock {\em Acad. Sci. USSR J. Phys.}, 10:13--24, 1946.
\newblock (Engl. transl).

\bibitem{Levy2000}
M.~Levy.
\newblock {\em Parabolic Equation Methods for Electromagnetic Wave
  Propagation}.
\newblock Electromagnetic Waves. Institution of Engineering and Technology,
  2000.

\bibitem{Lindman1975}
E.~L. Lindman.
\newblock ''{F}ree-space'' boundary vonditions for the time dependent wave
  equation.
\newblock {\em Journal of Computational Physics}, 18(1):66 -- 78, 1975.

\bibitem{Engquist1977}
B.~Engquist and A.~Majda.
\newblock Absorbing boundary conditions for the numerical simulation of waves.
\newblock {\em Math. Comp.}, 31:629--651, 1977.

\bibitem{Ristow1994}
D.~Ristow and T.~Rühl.
\newblock Fourier finite‐difference migration.
\newblock {\em Geophysics}, 59(12):1882--1893, 1994.

\bibitem{Li1991}
Z.~Li.
\newblock Compensating finite-difference errors in 3-d migration and modeling.
\newblock {\em Geophysics}, 56(10):1650--1660, 1991.

\bibitem{Zhang2008}
J.-H. Zhang, W.-M. Wang, L.-Y. Fu, and Z.-X. Yao.
\newblock 3{D} {F}ourier finite-difference migration by
  alternating-direction-implicit plus interpolation.
\newblock {\em Geophysical Prospecting}, 56(1):95--103, 2008.

\bibitem{ZHAO2019118}
H.~Zhao, L.-J. Gelius, M.~Tygel, E.~{Harris Nilsen}, and A.~{Kjelsrud Evensen}.
\newblock 3d prestack fourier mixed-domain (fmd) depth migration for vti media
  with large lateral contrasts.
\newblock {\em Journal of Applied Geophysics}, 168:118--127, 2019.

\bibitem{Gazdag1984}
J.~Gazdag and P.~Sguazzero.
\newblock Migration of seismic data by phase shift plus interpolation.
\newblock {\em Geophysics}, 49(2):124--131, 1984.

\bibitem{Stoffa1990}
P.~L. Stoffa, J.~T. Fokkema, R.~M. de~Luna~Freire, and W.~P. Kessinger.
\newblock Split‐step fourier migration.
\newblock {\em Geophysics}, 55(4):410--421, 1990.

\bibitem{Margrave1999}
G.~F. Margrave and R.~J. Ferguson.
\newblock Wavefield extrapolation by nonstationary phase shift.
\newblock {\em Geophysics}, 64(4):1067--1078, 1999.

\bibitem{Zhang2005}
Y.~Zhang, G.~Zhang, and N.~Bleistein.
\newblock Theory of true-amplitude one-way wave equations and true-amplitude
  common-shot migration.
\newblock {\em Geophysics}, 70(4):E1--E10, 2005.

\bibitem{Lee01101985}
M.~W. Lee and S.~Y. Suh.
\newblock Optimization of one-way wave equations.
\newblock {\em Geophysics}, 50(10):1634--1637, 1985.

\bibitem{Milinazzo1997}
F.~A. Milinazzo, C.~A. Zala, and G.~H. Brooke.
\newblock Rational square-root approximations for parabolic equation
  algorithms.
\newblock {\em The Journal of the Acoustical Society of America},
  101(2):760--766, 1997.

\bibitem{Lu1998}
Y.~Y. Lu.
\newblock A complex coefficient rational approximation of 1+x.
\newblock {\em Applied Numerical Mathematics}, 27(2):141--154, 1998.

\bibitem{Yevick2000}
D.~Yevick and D.~J. Thomson.
\newblock Complex {P}ad\'{e} approximants for wide-angle acoustic propagators.
\newblock {\em The Journal of the Acoustical Society of America},
  108(6):2784--2790, 2000.

\bibitem{Amazonas2010}
D.~Amazonas, R.~Aleixo, J.~Schleicher, and J.~Costa.
\newblock Anisotropic complex {P}ade hybrid finite-difference depth migration.
\newblock {\em Geophysics}, 75(2):S51, 2010.

\bibitem{Marchuk1968}
G.~I. Marchuk.
\newblock Some application of splitting-up methods to the solution of
  mathematical physics problems.
\newblock {\em Applik Mat}, 13(2):103--132, 1968.

\bibitem{Strang1968}
G.~Strang.
\newblock On the construction and comparison of difference schemes.
\newblock {\em SIAM Journal on Numerical Analysis}, 5(3):506--517, 1968.

\bibitem{Collino1995}
F.~Collino and P.~Joly.
\newblock Splitting of operators, alternate directions, and paraxial
  approximations for the three-dimensional wave equation.
\newblock {\em SIAM Journal on Scientific Computing}, 16(5):1019--1048, 1995.

\bibitem{Costa2013}
C.~A.~N. Costa, I.~S. Campos, J.~C. Costa, F.~A. Neto, J.~Schleicher, and
  A.~Novais.
\newblock Iterative methods for 3d implicit finite-difference migration using
  the complex pade approximation.
\newblock {\em Journal of Geophysics and Engineering}, 10(4):045011, 2013.

\bibitem{Terekhov2017}
A.~V. Terekhov.
\newblock The {L}aguerre finite difference one-way equation solver.
\newblock {\em Computer Physics Communications}, 214:71 -- 82, 2017.

\bibitem{Terekhov2018}
A.~V. Terekhov.
\newblock The stabilization of high-order multistep schemes for the {L}aguerre
  one-way wave equation solver.
\newblock {\em Journal of Computational Physics}, 368:115--130, 2018.

\bibitem{Bamberger1988a}
A.~Bamberger, B.~Engquist, L.~Halpern, and P.~Joly.
\newblock Parabolic wave equation approximations in heterogenous media.
\newblock {\em SIAM Journal on Applied Mathematics}, 48(1):99--128, 1988.

\bibitem{Temme1990}
N.~M. Temme.
\newblock Asymptotic estimates for {Laguerre} polynomials.
\newblock {\em Zeitschrift f{\"u}r angewandte Mathematik und Physik},
  41(1):114--126, Jan 1990.

\bibitem{Ernst2012}
O.~G. Ernst and M.~J. Gander.
\newblock {\em Why it is Difficult to Solve Helmholtz Problems with Classical
  Iterative Methods}, pages 325--363.
\newblock Springer Berlin Heidelberg, Berlin, Heidelberg, 2012.

\bibitem{Integral_Transform}
L.~Debnath and D.~Bhatta.
\newblock {\em Integral Transforms and Their Applications, Second Edition}.
\newblock Taylor \& Francis, 2006.

\bibitem{Terekhov2019}
A.~V. Terekhov.
\newblock Generating laguerre expansion coefficients by solving a
  one-dimensional transport equation (https://arxiv.org/abs/1809.06794), 2018.

\bibitem{Sumita1988}
U.~Sumita and M.~Kijima.
\newblock Theory and algorithms of the {L}aguerre transform, part {I}:{T}heory.
\newblock {\em Journal of the Operations Research Society of Japan},
  31(4):467--495, 1988.

\bibitem{Nussbaumer1982}
H.~J. Nussbaumer.
\newblock {\em Fast Fourier Transform and Convolution Algorithms}.
\newblock Springer-Verlag, 1982.

\bibitem{Tam1993262}
C.~K.~W. Tam and J.~C. Webb.
\newblock Dispersion-relation-preserving finite difference schemes for
  computational acoustics.
\newblock {\em Journal of Computational Physics}, 107(2):262--281, 1993.

\bibitem{Zhang2013511}
J.-H. Zhang and Z.-X. Yao.
\newblock Optimized explicit finite-difference schemes for spatial derivatives
  using maximum norm.
\newblock {\em Journal of Computational Physics}, 250:511--526, 2013.

\bibitem{CrankNicolson2}
J.~Crank and P.~Nicolson.
\newblock A practical method for numerical evaluation of solutions of partial
  differential equations of the heat-conduction type.
\newblock {\em Proc. Camb. Phil. Soc.}, 43(1):50--67, 1947.

\bibitem{Golub1989history}
G.~H. Golub and D.~P. O’Leary.
\newblock Some history of the conjugate gradient and {L}anczos algorithms:
  1948–1976.
\newblock {\em SIAM Review}, 31(1):50--102, 1989.

\bibitem{VanDerVorst_2003}
H.~A. van~der Vorst.
\newblock {\em Iterative Krylov Methods for Large Linear Systems}, volume~13.
\newblock Cambridge University Press, 2003.

\bibitem{Davis2007}
P.~J. Davis and P.~Rabinowitz.
\newblock {\em Methods of Numerical Integration}.
\newblock Dover Books on Mathematics Series. Dover Publications, 2007.

\bibitem{Terekhov2020}
A.~V. Terekhov.
\newblock An extra-components method for evaluating fast matrix-vector
  multiplication with special functions (https://arxiv.org/abs/2004.11610),
  2020.

\bibitem{claerbout1976}
J.~F. Claerbout.
\newblock {\em Fundamentals of geophysical data processing with applications to
  petroleum prospecting}.
\newblock McGraw-Hill, New York, March 1976.

\bibitem{Bunks1995}
C.~Bunks.
\newblock Effective filtering of artifacts for implicit finite-difference
  paraxial wave equation migration.
\newblock {\em Geophysical Prospecting}, 43(2):203--220, 1995.

\bibitem{Zhou1999}
H.~Zhou and G.~A. McMechan.
\newblock Parallel {B}utterworth and {C}hebyshev dip filters with applications
  to 3-{D} seismic migration.
\newblock {\em Geophysics}, 64(5):1573--1578, 1999.

\bibitem{Paffenholz}
J.~Paffenholz, B.~McLain, J.~Zaske, and P.~J. Keliher.
\newblock {\em Subsalt multiple attenuation and imaging: Observations from the
  Sigsbee2B synthetic dataset}, chapter 538, pages 2122--2125.

\bibitem{Gray2001}
S.~H. Gray, J.~Etgen, J.~Dellinger, and D.~Whitmore.
\newblock Seismic migration problems and solutions.
\newblock {\em Geophysics}, 66(5):1622--1640, 2001.

\bibitem{Valenciano2005}
A.~A. Valenciano and B.~Biondi.
\newblock {\em 2‐D deconvolution imaging condition for shot‐profile
  migration}, pages 1059--1062.
\newblock 2005.

\bibitem{Guitton2007}
A.~Guitton, A.~Valenciano, D.~Bevc, and J.~Claerbout.
\newblock Smoothing imaging condition for shot-profile migration.
\newblock {\em Geophysics}, 72(3):S149--S154, 2007.

\bibitem{Schleicher2008}
J.~Schleicher, J.~C. Costa, and A.~Novais.
\newblock A comparison of imaging conditions for wave-equation shot-profile
  migration.
\newblock {\em Geophysics}, 73(6):S219--S227, 2008.

\bibitem{GaussianBeams_popov}
M.~M. Popov.
\newblock A new method of computation of wave fields using gaussian beams.
\newblock {\em Wave Motion}, 4, 1982.

\bibitem{Cerveny1985}
V.~Cerveny.
\newblock Gaussian beam synthetic seismograms.
\newblock {\em J. Geophys.}, 58:44--72, 1985.

\bibitem{Golub1989}
G.H. Golub and C.F. Van~Loan.
\newblock {\em Matrix computations (3rd ed.)}.
\newblock Johns Hopkins University Press, Baltimore, MD, USA, 1996.

\bibitem{lahaye_modern_2017}
D.~Lahaye, J.~Tang, and K.~Vuik, editors.
\newblock {\em Modern {Solvers} for {Helmholtz} {Problems}}.
\newblock Geosystems {Mathematics}. Springer International Publishing, Cham,
  2017.

\bibitem{aminzadeh19973}
F.~Aminzadeh, B.~Jean, and T.~Kunz.
\newblock {\em 3-{D} salt and overthrust models}.
\newblock Society of Exploration Geophysicists, 1997.

\bibitem{Cerjan1985}
C.~Cerjan, D.~Kosloff, R.~Kosloff, and M.~Reshef.
\newblock A nonreflecting boundary condition for discrete acoustic and elastic
  wave equations.
\newblock {\em Geophysics}, 50(4):705--708, 1985.

\bibitem{Terekhov2015206}
A.~V. Terekhov.
\newblock Spectral-difference parallel algorithm for the seismic forward
  modeling in the presence of complex topography.
\newblock {\em Journal of Applied Geophysics}, 115(0):206--219, 2015.

\end{thebibliography}
\end{document}